\documentclass[12pt,psamsfonts]{amsart}
\usepackage{amsthm,amsmath,amsfonts,amssymb,amsthm,graphicx,epsf,epstopdf,color,pifont,subfigure,flafter,bm,amscd}
\usepackage{graphicx,overpic}
\usepackage[margin=3.6cm]{geometry}
\usepackage{ifpdf}

\newtheorem {theorem} {Theorem} %[section]
\newtheorem {proposition} [theorem] {Proposition}

\newtheorem {lemma} [theorem] {Lemma}

\newcommand{\R}{\ensuremath{\mathbb{R}}}

\begin{document}

\title[A Class of Quintic Polynomial Differential System]
{Global Centers in a Class of Quintic Polynomial Differential System}

\author[L.P.C. da Cruz] {Leonardo P.C. da Cruz$^1$}
\address{$^1$ Departamento de Matem\'atica, Instituto de Ci\^encias Matem\'aticas e Computa\c c\~ao, Universidade de S\~ao Paulo, Avenida Trabalhador S\~ao Carlense, 400, 13566-590, S\~ao Carlos, SP, Brazil}
\email{leonardocruz@icmc.usp.br, regilene@icmc.usp.br}

\author[J. Llibre] {Jaume LLibre $^2$}
\address{$^2$ Departament de Matem\`{a}tiques, Universitat Aut\`{o}noma de Barcelona, 08193 Bellaterra, Barcelona, Catalonia, Spain}
\email{jaume.llibre@mat.uab.cat}



\begin{abstract}
A center of a differential system in the plane $\R^2$ is an equilibrium point $p$ having a neighborhood $U$ such that $U\setminus \{p\}$ is filled of periodic orbits. A center $p$ is global when $\R^2\setminus \{p\}$ is filled of periodic orbits. In general is a difficult problem to distinguish the centers from the foci for a given class of differential systems, and also it is difficult to distinguish the global centers inside the centers.
The goal of this paper is to classify the centers and the global centers of the following class of quintic polynomial differential systems
$$
\dot{x}= y,\quad \dot{y}=-x+a_{05}\,y^5+a_{14}\,x\,y^4+a_{23}\,x^2\,y^3+a_{32}\,x^3\,y^2+a_{41}\,x^4\,y+a_{50}\,x^5,
$$
in the plane $\R^2$.
\end{abstract}

\maketitle
	
\section{Introduction and statement of the main results}

We consider polynomial differential systems
\begin{equation}\label{eq:-1}
(\dot{x},\dot{y})=(P(x,y),Q(x,y)),
\end{equation}
defined in the plane $\R^2$. Here the dot denotes derivative with respect to the time $t$. We are interested in the subclass of these polynomial differential systems having an equilibrium point whose linear part has eigenvalues purely imaginary. After an affine change of variables and a rescaling of the time (if necessary) such polynomial differential systems can be written into the form
 \begin{equation}\label{eq:0}
(\dot{x},\dot{y})=(-y+P_n(x,y),x+Q_n(x,y)),
\end{equation}
where $P_n$ and $Q_n$ are polynomials of degree $n$, which do not have neither constant nor linear terms.

The problem of distinguishing whether the equilibrium point at the origin of system \eqref{eq:0} is a center or a focus is a classical problem, known as the  \emph{center-focus problem}. Even this problem was partially solved by Lyapunov, see \cite{Lia1893}, it has been studied for some fixed values of the degree $n$ during more than a century by many authors. The only family completely investigated is the family of the polynomial differential systems of degree $2$, denoted simply by quadratic systems. The study of this family was started by Dulac in 1908 in \cite{Dul1908}, and also performed by Kapteyn some years later, see \cite{Kap1911,Kap1912}. Up to the work of Frommer \cite{Frommer1934} the conditions for the existence of a center in that family were not published. The correct center conditions were published by Saharnikov \cite{Saha1948} and later by Sibirski{\u\i} \cite{Sibi1954,Sibi1955}. The center conditions are simpler and the center-focus problem is easier to solve if the system is written in complex coordinates, see \cite{Zol1994}. For the complete cubic family (when in systems \eqref{eq:0} we have $n=3$), the problem remains unsolved.

In 1992 Galeotti and Villarini in \cite{MGaMVill1992}, proved that if the degree $n$ of a polynomial differential system \eqref{eq:0} is even such a system cannot have a global center. Recently Llibre and  Valls proved the same in an easier way, see \cite{LliVall2021}. In \cite{Conti1998} Conti proposed the following  problem: \emph{To classify all polynomial differential systems of degree odd having global centers}, this is a very difficult problem.

It is well known that all the centers of the polynomial differential systems \eqref{eq:0} with $n=1$ (i.e. of the linear differential systems) are global.

For $n=3$ the global linear centers (i.e. the centers with purely imaginary eigenvalues) and the global nilpotent centers (i.e. the centers having both eigenvalues zero but its linear part is not identically zero) having only homogeneous nonlinearities of degree $3$ are completely classified in \cite{JohJauValls2019} and \cite{JohJauValls2020}, respectively. Moreover also for $n=3$ in \cite{LuisJauValls2020} are classified the global centers of the Hamiltonian systems that are reversible with respect to the $x$-axis, and in \cite{LuisJauValls2022} are classified the global nilpotent centers of the cubic Hamiltonian systems.

When $n=5$ and for the systems \eqref{eq:0} with homogeneous nonlinearities of degree $5$ in \cite{LlibVall2023} are classified their global linear centers for the subclass of such systems that are reversible with respect to the $x$-axis, and in \cite{JohJauValls2023} with respect to the $y$-axis.

A natural continuation in the classification of the global centers is to consider another family of polynomial differential systems \eqref{eq:0} with $n=5$. Thus in this paper we consider the following class of quintic polynomials differential systems
\begin{equation}\label{eq:1}
\dot{x}= y,\quad
\dot{y}=-x+a_{05}\,y^5+a_{14}\,x\,y^4+a_{23}\,x^2\,y^3+a_{32}\,x^3\,y^2+a_{41}\,x^4\,y+a_{50}\,x^5,
\end{equation}
in the plane $\R^2$. Then the main result of this paper is the following one.

\begin{theorem}\label{thmmainn_1}
The polynomial differential systems \eqref{eq:1} have a global center at the origin of coordinates if, and only if, the following conditions hold $a_{32}\leq0, a_{50}\leq 0  \hspace{0.1cm} \text{and} \hspace{0.2cm} a_{14}=a_{41}=a_{23}=a_{05}=0.$
\end{theorem}

Before starting with the classification of the global center it is necessary to solve the center-focus problem. The classification of the centers has been possible despite the number six of parameters appearing. Thus in the following result we provide the answer to the center-focus classification problem.

\begin{theorem}\label{thmmainn_2}
The polynomial differential systems \eqref{eq:1} have a center at the origin of coordinates and the origin is the unique equilibrium point of the system if, and only if, the following conditions hold
$a_{50}\leq0 \hspace{0.1cm} \text{and} \hspace{0.2cm}  a_{41}=a_{23}=a_{05}=0.$	
\end{theorem}

This paper is structured as follows. In section \ref{se:preliminaries} we provide the necessary definitions, results and the algorithm to obtain the coefficients of the return map, the so-called Lyapunov constants. In sections \ref{se:necsuficondition} and \ref{se:bifurcation} we prove the Theorems \ref{thmmainn_2} and Theorem \ref{thmmainn_1}, respectively.

\section{Preliminaries}\label{se:preliminaries}
	
In this section we recall some classical concepts necessary to state and prove the results of this paper.

In subsection \ref{se:centerconditions} we remember how to obtain the coefficients of the return map (i.e. the Lyapunov constants) near a monodromic equilibrium point (i.e. a focus or a center). \\
In subsection \ref{se:bautincenter} we define the Bautin ideal and the center variety.\\
In subsection \ref{se:compactificacion} we recall the Poincar\'e compactification. \\
Finally in subsection \ref{se:globalcenter} we remember the result for classifying the global centers.

\subsection{The center conditions}\label{se:centerconditions}

We consider a polynomial differential system of degree $n$ with an equilibrium point at the origin of coordinates having its Jacobian matrix purely imaginary eigenvalues, i.e. we consider a system \eqref{eq:0}. So the origin is a focus or a center.

A non-constant analytical function defined in a neighborhood $\Omega$ of the origin, $H:\Omega\subset\mathbb{R}^2\rightarrow \mathbb{R}^2,$ is a first integral of system \eqref{eq:0} if it is constant along any solution $\gamma$ or, equivalently,
\begin{equation}\label{eq:3}
\left. \frac{\partial H}{\partial x} \dot{x}+\frac{\partial H}{\partial y}\dot{y}\right|_{\gamma}\equiv 0.
\end{equation}
In order to distinguish when the origin is a center we shall use the  Poincar\'e--Lyapunov Theorem, see \cite{IliYak2008, Liapounoff1965, Poincare1886, RomaSha2009}:

\begin{theorem}
The polynomial differential system \eqref{eq:0} has a center at the origin if and only if it admits a local analytic first integral of the form
\begin{equation}\label{eq:2_1}
H(x,y)=x^2+y^2 + \sum_{p=3}^\infty H_p(x,y), \; \mbox{where} \; H_p(x,y)=\sum_{\ell=0 }^p q_{p-\ell,\ell}x^{p-\ell}y^\ell.
\end{equation}
In addition, the existence of a formal first integral $H$ of the above form implies the existence of a local analytic first integral.
\end{theorem}

The necessary conditions for the existence of a first integral \eqref{eq:2_1} for system \eqref{eq:0} are obtained by looking for a formal series \eqref{eq:2_1} satisfying \eqref{eq:3}. Although \eqref{eq:3} is not always satisfied, it is always possible to choose coefficients of the formal power series \eqref{eq:2_1} that satisfy the following equation:
\begin{equation}\label{eqliapu}
\frac{\partial H}{\partial x} \dot{x}+\frac{\partial H}{\partial y}\dot{y}=\sum_{j=1}^\infty L_j (x^2+y^2)^{j+1},
\end{equation}
see \cite{BlowLlo1984,Lia1966}. We remark that any non-zero $L_j$ obstructs the origin to be a center. Then when at least one $L_j$ is different from zero it is a Lyapunov constant in a neighborhood of the origin. Hence system \eqref{eq:0} has no local analytic first integral and we say that the equilibrium point is a weak focus of order $k$ if the first non-zero coefficient in \eqref{eqliapu} is $L_k.$ The coefficient $L_j$ in \eqref{eqliapu} is called the $j$-th Lyapunov constant. The stability of the origin is given by the sign of the first non-zero $L_j.$ Moreover, note that the constants $L_j$ are rational functions whose numerators are polynomials depending on the coefficients of the polynomial system \eqref{eq:0}.

In order to compute the first $N$ Lyapunov constants we need to arrive to compute until the terms of order $2N+2$ in the series \eqref{eq:2_1}, i.e.
\begin{equation}\label{eq:4}
\widetilde{H}(x,y)=\dfrac{x^2+y^2}{2}+\sum_{p= 3}^{2N+2}\sum_{\ell=0 }^p q_{p-\ell,\ell}x^{p-\ell}y^\ell.
\end{equation}
After for each $i=3,\ldots,2N+2$ we equate to zero the coefficients of terms of degree $i$ in the expression
\begin{equation*}
\frac{\partial \widetilde H}{\partial x}\dot{x}+\frac{\partial \widetilde H}{\partial y}\dot{y}=\left(-y+\displaystyle\sum_{k=2}^n P_k(x,y)\right)\frac{\partial \widetilde H}{\partial x}+\left(x+\displaystyle\sum_{k=2}^n Q_k(x,y)\right)\frac{\partial \widetilde H}{\partial y}.
\end{equation*}
Hence starting with $i=3$ we solve in a recurrent way each linear system of $i+1$ equations with $i+1$ variables, $q_{p-\ell,\ell}$ for $\ell=0,\dots,p$. All linear systems corresponding to odd degrees, $i=2j+1,$ have a unique solution in terms of the previous values of $q_{p-\ell,\ell}$. As the determinant of the linear system that corresponds to an even degree, $i=2j+2,$ vanishes, we need to add an extra condition so that the linear system has a unique solution. In fact, at this step, we have one equation more than the number of variables. We add suitable equations, for the terms $(x^{2}+y^{2})^{j+2}$ for example, so that the derivative over the associated vector field becomes
\begin{equation}\label{eq:6}
\frac{\partial H}{\partial x}\dot{x}+\frac{\partial H}{\partial y}\dot{y}=\sum_{j=1}^\infty L_{j}(x^{2}+y^{2})^{j+2}.
\end{equation}
Therefore we define the Lyapunov constants associated with the extra conditions given above. In this context it is well-known that the first non-vanishing coefficient of \eqref{eq:6} has an odd subindex, and $L_{2k+1}$ is called the $k$th-order Lyapunov constant of system \eqref{eq:0}. An interesting property, described in \cite{Rou1998} and proved in \cite{CimGasMan2020}, of these Lyapunov constant is that for each $k$ we have that the ideals
\begin{equation}\label{cst}
\langle L_2,L_4,\ldots, L_{2j}\rangle \subset \langle L_3,L_5,\ldots, L_{2j-1}\rangle.
\end{equation}

\subsection{The Bautin ideal and the center variety}\label{se:bautincenter}

We called $\lambda \in \mathbb{R}^M$ the set of parameters of system \eqref{eq:0}, then the $L_j$'s are polynomials in $\lambda$. In addition, the set $\mathcal{B}^{\mathbb{R}}=\langle L_1,L_2,\ldots\rangle$ is an ideal in the polynomial ring $\mathbb{\mathbb{R}}[\lambda]$. The importance of the ideal $\mathcal{B}^{\mathbb{R}}$ follows from the fact that if all the $L_j$'s generating the ideal vanishes, then all the Lyapunov constants vanish and it is not necessary to compute all of them. Then we introduce the next definition which recalls the notion of \emph{Bautin ideal} and the \emph{center variety}.

The ideal defined by the Lyapunov constants $\mathcal{B}^{\mathbb{R}}=\langle L_1,L_2,\ldots\rangle\subset\mathbb{R}[\lambda]$ is called the {\it Bautin ideal}. The affine variety $\mathbf{V}^{\mathbb{R}}=\mathbf{V} (\mathcal{B}^\mathbb{R})$ is called the {\it center variety} of system \eqref{eq:0}, i.e.
$\mathbf{V}(\mathcal{B}^{\mathbb{R}})=\{\lambda \in \mathbb{R}^M: r(\lambda)=0, \; \forall \; r \in \mathcal{B}^{\mathbb{R}}\}$.

When we can explicitly determine the center variety we have the center-focus problem solved for system \eqref{eq:0}. However, in most of the cases, this is not a simple problem. On the other hand, the Hilbert Basis Theorem \cite{AdaLou1994,CoxLitOsh1997,MacBir1967} assures that $V(B^\R)$ is finitely generated. Then there exists a positive integer $j$ such that $\mathcal{B}^{\mathbb{R}}=\mathcal{B}^{\mathbb{R}}_j=\langle L_1,\ldots,L_j\rangle.$ In other words, we know that for $j$ big enough, the above algorithm provides a necessary set of conditions $\{L_j=0: j=1,\ldots, N\}$ in order that system \eqref{eq:0} be a center. The main difficulty follows from the fact that there is no technique to obtain $j$ a priori.

We can also say that the polynomials $L_j$ represent obstacles to the existence of a first integral. In particular, system~\eqref{eq:0} admits a first integral of the form~\eqref{eq:4} if and only if $L_j=0,$ for all $j\geq1.$ Thus the simultaneous vanishing of all focus quantities provides conditions that characterize when a system of the form~\eqref{eq:0} has a center at the origin. Note that the inclusion $\mathbf{V}^{\mathbb{R}}= \mathbf{V}(\mathcal{B}^\mathbb{R})\supset\mathbf{V}(\mathcal{B}^\mathbb{R}_j)$ holds for any $j\geq1.$ The opposite inclusion, for a fixed $j$, is verified by finding the irreducible decomposition of $\mathbf{V}(\mathcal{B}^ \mathbb{R}_j),$ see \cite{RomaSha2009}, then any point of each component of the decomposition corresponds to a system having a center at the origin.

\subsection{The Poincaré compactification}\label{se:compactificacion}

Given a planar polynomial differential system \eqref{eq:-1}, one crucial problem of the qualitative theory is to characterize the phase portraits in the Poincaré disc of this system. First we must characterize the local phase portraits of the finite and infinite equilibrium points to reach this challenging goal. In the theory of the Poincaré compactification the circle of the Poincaré disc represents the infinity of the plane $\R^2.$ See more details \cite{DumLliArt2006}.

For studying the equilibrium points on the circle at infinity, we need four local charts $U_1= \{(x,y): x>0\}$,  $V_1= \{(x,y): x<0\}$, $U_2= \{(x,y): y>0\}$, and $V_2= \{(x,y): y<0\}$.

To study the dynamics at infinity we need the followging expressions of the polynomial system \eqref{eq:-1} of degree $n$ on the Poincaré disc
\begin{equation}\label{eq:u1u2}
\begin{array}{l}
(\dot{x},\dot{y})= \Big(y^n\big(Q(1/y,x/y)-x\,P(1/y,x/y)\big), y^{n+1}P(1/y,x/y)\Big) \ \ \text{in} \ \ U_1,\\\\
(\dot{x},\dot{y})= \Big(y^n\big(P(1/y,x/y)-x\,Q(1/y,x/y)\big),y^{n+1} Q(1/y,x/y)\Big)\ \ \text{in} \ \ U_2.
\end{array}
\end{equation}
The expressions of system \eqref{eq:-1} in the local charts $V_i$ for $i=1,2$ are the same than in the local charts $U_i$ for $i=1,2$ multiplied by $(-1)^n$.

For studying the infinite equilibrium points which in all these local charts are of the form $(x,0)$ it is sufficient to study the infinite equilibrium points of the local chart $U_1$ and the origin of the local chart $U_2$, again for more details see Chapter 5 of \cite{DumLliArt2006}.

\subsection{The classification of the global centers}\label{se:globalcenter}

The following result gives the conditions in order that a polynomial differential system in the plan $\R^2$ has a global center, see \cite{LlibVall2023}.

\begin{proposition}\label{cg}	
A polynomial differential system of degree $n$ in $\R^2$ without a line of equilibrium points at infinity, has a global center if, and only if, it has a unique finite equilibrium point which is a center and all the local phase portraits of the infinite equilibrium points (if they exist) are formed by two hyperbolic sectors having all of them both separatrices on the infinite circle.	
\end{proposition}

\section{The classification of the centers}\label{se:necsuficondition}	

We devote this section to prove Theorem \ref{thmmainn_2}. As the proof is quite long, we have divided it in two propositions and a lemma. In the Proposition \ref{cent1} we prove the necessary conditions for having a center at the origin of coordinates, and in the Proposition~\ref{cent2} we establish sufficient conditions in order that the centers are global centers. Finally in Lemma \ref{l1} we give the conditions in order that the unique equilibrium point of the system be the origin.
	
\begin{proposition}\label{cent1}
If the origin of the quintic polynomial differential system \eqref{eq:1} is a center, then the parameters $a_{ij},$ with $i+j=5$ satisfy the conditions given in the statement of Theorem $\ref{thmmainn_2}$.
\end{proposition}

\begin{proof}
The trace and the determinant of the Jacobian matrix at the origin of system \eqref{eq:1} are zero and a positive, respectively. So the origin is a focus or a center. We need to compute some Lyapunov constants for distinguish the centers from the foci, and since the system has six parameters we must compute at least six Lyapunov constants. Then we have the following system of equations with six parameters
\begin{equation*}\label{S}
\mathcal{S}=\{L_1=L_2=\ldots=L_{2k}=L_{2k+1}=0\}, \hspace{0.1cm} \text{and} \hspace{0.2cm} k\geq6.
\end{equation*}

Following the approach described in subsection \ref{se:centerconditions} for the computation of the center conditions $L_i$, we have computed the Lyapunov constants $L_i$ for $i=1,\ldots,17$, and according with the property \eqref{cst} we must solve the following algebraic system composed only with odd Lyapunov constants
\begin{equation*}\label{P}
\mathcal{Q}=\{L_1=L_3=\ldots=L_{2k+1}=0\}, \hspace{0.1cm} \text{and} \hspace{0.2cm} k=8,
\end{equation*}
where $L_1\equiv0,$ and this Lyapunov constants are polynomial in the parameters $a_{ij},$ with $i+j=5.$ Due to the huge expressions of these Lyapunov constants, we only provide in what follows the first four Lyapunov constants, where we have denoted $a_{ij}$ by $a_{i}$ for $i=0,\ldots,5$:
\begin{equation*}
\begin{aligned}
L_3=&-\frac{1}{16}\big(5a_{0}+a_{4}+a_{2}\big),\\
L_5=&-\frac{1}{128}\big(5a_{5}a_{4} + \frac{7}{5}a_{4}a_{3}+ \frac{1}{5}a_{4}a_{1}+ \frac{3}{2}a_{5}a_{2}+\frac{9}{10}a_{3}a_{2}+\frac{7}{10}a_{2}a_{1}\big),\\
L_7=&-\frac{1}{2400}\big(a_{4}^3-\frac{9}{32}a_{2}^3-\frac{ 765}{784}a_{3}^2a_{2}+\frac{165275}{3136}a_{5}^2a_{4} + \frac{825}{196}a_{5}^2a_{2} + \frac{11625}{1568}a_{5}a_{4}a_{3}\big.\\
& \big.+ \frac{21075}{1568}a_{5}a_{4}a_{1} + \frac{375}{196}a_{5}a_{3}a_{2}  - \frac{5}{8}a_{4}^2a_{2} - \frac{155}{448}a_{4}a_{3}^2 +\frac{ 3825}{1568}a_{4}a_{3}a_{1} - \frac{9}{8}a_{4}a_{2}^2\\
& - \frac{675}{448}a_{4}a_{1}^2 \big),\\
L_9=&\phantom{-}\frac{1}{36864}\big(a_{3}^3a_{2}-\frac{27665}{27}a_{5}^3a_{4}-\frac{440}{9}a_{5}^3a_{2}-\frac{7493}{21}a_{5}^2a_{4}a_{3}-\frac{ 2707}{9}a_{5}^2a_{4}a_{1}\big.\\ 
&+\frac{404}{63}a_{5}^2a_{3}a_{2}-\frac{1363574}{46305}a_{5}a_{4}^3 + \frac{388547}{46305}a_{5}a_{4}^2a_{2}- \frac{4813}{108}a_{5}a_{4}a_{3}^2 - \frac{73}{36}a_{4}a_{3}^3\\
&- \frac{1700}{21}a_{5}a_{4}a_{3}a_{1}+ \frac{16031}{1470}a_{5}a_{4}a_{2}^2 + \frac{47}{7}a_{5}a_{3}^2a_{2}+\frac{517}{420}a_{5}a_{2}^3-\frac{ 53482}{11025}a_{4}^3a_{3}\\
&\big.-\frac{1388}{343}a_{4}^3a_{1}+ \frac{92923}{77175}a_{4}^2a_{3}a_{2}
- \frac{85}{12}a_{4}a_{3}^2a_{1} + \frac{3907}{2450}a_{4}a_{3}a_{2}^2  + \frac{209}{700}a_{3}a_{2}^3\big),\ldots.
\end{aligned}
\end{equation*}

Here we get $L_{2k+1}$ assuming that $L_{2k+1}\in\langle L_3\dots,L_{2k-1} \rangle$ for $k=1,\ldots,8$ and $L_{11}, L_{13},$ $L_{15}$ are polynomials of degrees 5, 6 and 7 in the variables $a_0,a_1,\ldots,a_5,$ respectively. Moreover, we get that $L_{17}\equiv0.$  Now we need to solve the algebraic system of equations $\mathcal{Q}$. However despite this system has only six variables and seven equations the usual mechanisms for solving it fail. Then we determine the irreducible components of the variety
\begin{equation}\label{V}
\mathbf{V}=\mathbf{V}(L_3,L_5,L_7,L_9,L_{11},L_{13},L_{15}),
\end{equation}
consists of using the Gianni–Trager–Zacharias algorithm, see \cite{GiaTraZac1988}, for determining the irreducible components of the variety \eqref{V}. The main function used is minAssGTZ, it is implemented in the library primdec.lib included on the algebraic computational system SINGULAR, see \cite{DecGrePfiSch2018,DecPfiSchLap2018}. Here using this algorithm we are able to compute the decomposition, and obtain the necessary conditions to have a center finding the irreducible decomposition of the variety. Working in $\mathbb{Q}[a_0,a_1,\ldots,a_5]$ the minimal correponding prime ideal of $\mathcal{R}=\langle L_3,L_5,\dots,L_{2k+1}\rangle$ with $k=7$ provided by SINGULAR is
\begin{equation}\label{Ts}
\begin{aligned}
\mathcal{T}_1=& \langle a_{0},a_{2},a_{4}+a_{2}+5a_{0}\rangle,\\	
\mathcal{T}_2=&\langle 18a_{3}^2+49a_{2}^2,a_{0}, a_{1},a_{4}+a_{2}+5a_{0},7a_{5}+a_{3}\rangle, \vspace{0.5cm} \\
\mathcal{T}_3=&\langle a_{2}^2+16a_{5}^2-20a_{2}a_{0}+100a_{0}^2, a_{3}a_{5}+6a_{5}^2-5a_{2}a_{0}+50a_{0}^2, a_{4}+a_{2}+5a_{0}+ \\
&  a_{3}a_{2}+6a_{2}a_{5}-10a_{3}a_{0}+20a_{5}a_{0} ,a_{3}^2-36a_{5}^2+60a_{2}a_{0}-200a_{0}^2, 5a_{5}+a_{3}+a_{5} \rangle.\vspace{0.5cm}\\
\end{aligned}
\end{equation}
The next step is to show that $\sqrt{\mathcal{T}}=\sqrt{\mathcal{R}},$ where
$\mathcal{T}=\bigcap_{k=1}^{3}\mathcal{T}_{k}$ in $\mathbb{Q}[a_0,a_1,\ldots,a_5].$ We denote by $\sqrt{U}$ the radical of the ideal $U$. In general, it is simpler to verify the double inclusion instead of computing the radicals. Adding a new artificial parameter $w$, this property can be seen checking that $\{1\}$ is the Gröbner basis of the next list of ideals, $\langle1-w L_{2k+1},\mathcal{T}\rangle,$ for $k=1,\ldots,7$ and $\langle1-w p,\mathcal{R}\rangle,$ for every $p \in \mathcal{T}$.

Finally we study the variety of each minimal prime ideal of \eqref{Ts}. Fristly, for $\mathcal{T}_1:$ the variety is given by the solution of the algebraic system $\{a_{0}=a_{2}=a_{4}+a_{2}+5\,a_{0}=0\},
$
then we get $a_{0}=a_{2}=a_{4}=0.$ Seconfly, for $\mathcal{T}_2:$ we must study the algebraic system 
$\{18\,a_{3}^2+49\,a_{2}^2=a_{0}=a_{1}=a_{4}+a_{2}+5\,a_{0}=7\,a_{5}+a_{3}=0\},$
and directaly, the variety belongs to the complex space since that the solution is complexa. Similarly, for $\mathcal{T}_3:$ we get that the variety is complex. Thus filtering these solutions we obtain the only  center condition given in the statement of Theorem \ref{thmmainn_2}, is given by variety of the minimal corresponding prime ideal $\mathcal{T}_1$. 

\end{proof} 

\begin{proposition}\label{cent2}
Under the conditions of Theorem~\ref{thmmainn_2} the quintic polynomial differential system \eqref{eq:1} has a center at the origin.
\end{proposition}

\begin{proof}
Since the origin of system \eqref{eq:1} is a focus or a center, and this system is invariant under the symmetry $(x,y,t)\rightarrow(x,-y,-t),$ it follows that the origin is a center.
\end{proof}
	
\begin{lemma}\label{l1}
Under the conditions $a_{41}=a_{23}=a_{05}=0$ the origen is the unique finite equilibrium point of the quintic polynomial differential system \eqref{eq:1}  if, and only if $a_{50}\leq0.$
\end{lemma}	

\begin{proof}
Assuming that $p=(\alpha,\beta)$ is an equilibrium point, we have $\beta=0$. So $p=(\alpha,0)$ and must satisfy the conditions $\alpha(a_{50}\alpha^4 - 1)=0.$ Therefore, $\alpha=0$ is the unique solution, if $a_{50}\leq0.$
\end{proof}	
	
\begin{proof}[Proof of Theorem $\ref{thmmainn_2}$]
It follows directly from Propositions~[\ref{cent1}, \ref{cent2}], and Lemma~\ref{l1}.
\end{proof}

\section{The global centers}\label{se:bifurcation}	
	
This section is devoted to prove Theorem \ref{thmmainn_1}. The proof of theorem follows using Proposition~\ref{cg}, i.e, we will give sufficient conditions to classify the global center. First we assume the conditions of the statement of Theorem~\ref{thmmainn_2}, so the origin of system \eqref{eq:1} is the unique equilibrium point and it is a center. So, we get the differential system
\begin{equation}\label{center}
	\dot{x}=y,\quad
	\dot{y}=-x+a_{14}\,x\,y^4+a_{32}\,x^3\,y^2-a^2\,x^5,
\end{equation}	
where $-a^2=a_{50}\leq0.$ So, from \eqref{eq:u1u2} and Theorem \ref{thmmainn_2}, the differential system \eqref{center} in the local charts $U_1$ and $U_2$ becomes
\begin{equation}\label{U1}
\dot{x}=-\,y^4\,x^2+a_{14}\,x^4-y^4+a_{32}\,x^2-a^2 ,\quad
\dot{y}=-\,y^5\,x,
\end{equation}	
and
\begin{equation}\label{U2}
\dot{x}= a^2\,x^6+ y^4\,x^2-a_{32}\,x^4+ y^4-a_{14}\,x^2 ,\quad
\dot{y}= x\,y\,(a^2\,x^4 + y^4-a_{32}\,x^2 - a_{14}),
\end{equation}	
respectively.

Before proving Theorem \ref{thmmainn_1} we define the concept of a {\it characteristic direction} at an equilibrium point and how to compuet them.
Given a polynomial differential system of degree $n,$ of the form
\begin{equation}\label{eq_dc}
(\dot{x},\dot{y})=\left(P_k(x,y)+ h.o.t., Q_k(x,y)+ h.o.t.\right),
\end{equation}
where $P_k(x,y)$ and $Q_k(x,y)$ are the terms of lower degree $k\ge 1,$ of the differential system \eqref{eq:-1}. Here $n\geq k,$ and h.o.t. denotes higher order terms. If the origin $p=(0,0)$ is an equilibrium point of system \eqref{eq_dc}, then the {\it characteristic direction} of the orbit $\gamma(t)$ at $p$ tending to $p$ in positive time (respectively in negative time) is the limit $\lim_{t\to \infty}(\gamma(t)-p)/\|\gamma(t)-p\|$ (respectively $\lim_{t\to-\infty}(\gamma(t)-p)/\|\gamma(t)-p\|$), if such a limit exists. Moreover, consider the homogeneous polynomial
\begin{equation}\label{cd}
\gamma_{k}= P_k(x,y)\,y-Q_k(x,y)\,x,
\end{equation}
the possible characteristic directions of the orbits starting or ending at the equilibrium point localized at the origin of coordinates are given by the real linear factors of the homogeneous polynomial \eqref{cd}. For more details on the characteristic directions, see \cite{AndLoeGorMai1973}.

\begin{proof}[Proof of Theorem $\ref{thmmainn_1}$]
We shall use Proposition~\ref{cg} in this proof, that gives the necessary and sufficient conditions for classifying global centers.
Thus we must determine the local phase portraits of the infinie equilibrium points of system \eqref{center}. 
%Here, in order to make it more short, in same time of this prove we call the origin of the coordinate space $(0,0)$ by $\mathcal{P}.$ 
The linear part of the infinite equilibrium point (the origin) in the chart $U_2$ of the system \eqref{U2} is identically zero. Thus in order to determine its local phase portrait we must do blow-ups.

Assume that $a_{14}\ne 0,$ in the system \eqref{U2}. From \eqref{cd} the characteristic directions at the origin are obtained from $\gamma_{2}=0,$ where $P_2(x,y)=-\,a_{14}\,x^2$ and $Q_2(x,y)=-\,a_{14}\,x\,y.$ So all directions are characteristic. Then we do the vertical blow-up $(x,y)\rightarrow (x_1,y_1\,x_1),$ to system \eqref{U2}, and we get
\begin{equation}\label{U2_1}
\dot{x_1}= x_1^2\,(x_1^4\,y_1^4+a^2\,x_1^4+x_1^2\,y_1^4-a_{32}\,x_1^2 -a_{14}), \quad \dot{y_1}=-x_1^3\,y_1^5,
\end{equation}
and rescaling the time $(x_2,y_2,t)\rightarrow(x_1,y_1,t/x_1^2)$ in \eqref{U2_1} we obtain
\begin{equation}\label{U2_2}
\dot{x_2}= x_2^4\,y_2^4+a^2\,x_2^4+x_2^2\,y_2^4-a_{32}\,x_2^2-a_{14}, \quad
\dot{y_2}=-x_2\,y_2^5.
\end{equation}
Going back through the changes of variables the local phase portrait at the equilibrium point (the origin) is shown in Figure~\ref{fi:star}, if $a_{14}>0$. When $a_{14}<0,$ the local phase portrait of the origin is the one of Figure~\ref{fi:star}, reversing the orientation of the orbits. So when $a_{14}\ne 0$ there are orbits which go or come from the infinity in system \eqref{center}, and consequently the center of this system cannot be global.

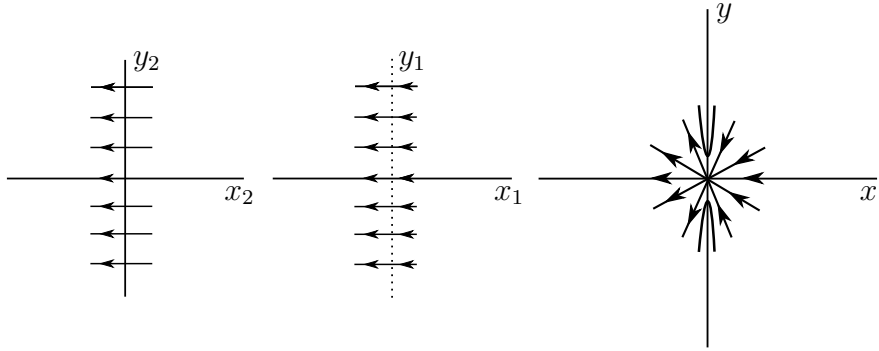
\begin{figure}[h]
\begin{center}
\begin{overpic}[height=4.5cm]{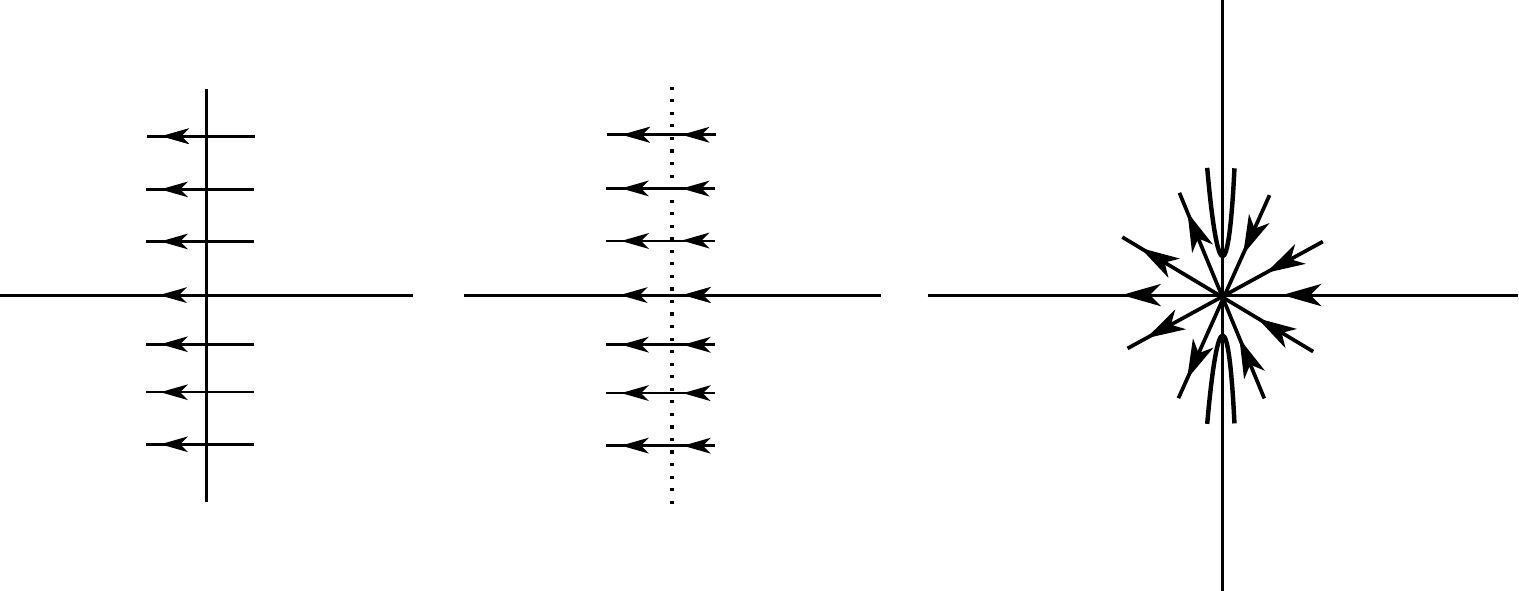}
\put(25,17){$x_2$}
\put(14.5,32.5){$y_2$}
\put(56,17){$x_1$}
\put(45,32.5){$y_1$}
\put(98,17){$x$}
\put(81.5,38){$y$}
\end{overpic}
\end{center}
\caption{Here $a_{14}>0$. In the picture on the left of the figure there is the local phase portrait in a neighborhood of the $y_2$-axis of systems \eqref{U2_2}. In the picture in the middle of the figure there is the local phase portrait in a neighborhood of the $y_2$-axis of system \eqref{U2_1}. The local phase portrait of the equilibrium point localized at the origin of the chart $U_2$ is shown in the rigth picture.} \label{fi:star}
\end{figure}
Assume now that $a_{14}=0$. Then system \eqref{U1} writes
\begin{equation}\label{U1_1}
\dot{x}=-\,y^4\,x^2-y^4+a_{32}\,x^2-a^2 ,\quad \dot{y}=-\,y^5\,x.
\end{equation}
The infinite equilibrium points of this system are
\begin{equation}\label{eq_U1_1}
\mathcal{P}_\pm=\left(\pm \frac{a}{\sqrt{a_{32}}},0\right),
\end{equation}
if they exist. When they exist, since  the differential system \eqref{U1_1} is invariant under the symmetry $(x,y,t)\to (-x,y,-t)$ we only need to study the local phase portrait at the infinite equilibrium point $\mathcal{P}_+$.

We divide the study of the possible infinite singular points of system \eqref{U1_1} into the following six cases:
\begin{equation*}\label{ci}
\begin{aligned}
c_1&=\{a_{32}>0,a\neq0\}, \ \ \ \  c_2=\{a_{32}>0,a=0\},\\
c_3&=\{a_{32}<0,a\neq0\}, \ \ \ \  c_4=\{a_{32}<0,a=0\},\\
c_5&=\{a_{32}=0,a\neq0\}, \ \ \ \  c_6=\{a_{32}=a=0\}.
 \end{aligned}
\end{equation*}	

$\bullet$ Case $c_1$.  Then translating the equilibrium $\mathcal{P}_+$ at the origin, the system \eqref{U1_1} becomes
$\dot{x}=2 a \sqrt{a_{32}}\, x+\ldots$, $\dot{y}=\ldots$,
here the dots mean terms of degree larger than one in the variables $x$ and $y.$ Then this equilibrium, by Theorem 2.15 of \cite{DumLliArt2006}, is a semi-hyperbolic saddle, or node, or  saddle-node, and consequently some orbit of the system \eqref{center} goes or comes from the infinity, and system \eqref{center} cannot have a global center.

$\bullet$ Case $c_2$. Then system \eqref{U2} writes
\begin{equation}\label{c2}
\dot{x}=x^2\,y^4+y^4-a_{32}\,x^4,\quad 
\dot{y}=y\,x\,(x^4-a_{32}\,x^2).
\end{equation}	
From \eqref{cd} the characteristic direction at origin is $\gamma_{4}=y^5=0$, and we can do the vertical blow-up $(x,y)\rightarrow(x_1,y_1\,x_1),$ without loosing information because $x=0$ is not a characteristic direction. So system \eqref{c2} goes over the system
\begin{equation}\label{c2_1}
\dot{x_1}=-x_1^4\,(-x_1^2\,y_1^4-y_1^4+a_{32}),\quad
\dot{y_1}=-x_1^3\,y_1^5,
\end{equation}	
and doing the rescaling of the time $(x_2,y_2,t)\to (x_1,y_1,t/x_1^3)$ in system \eqref{c2_1}, we get the system
\begin{equation}\label{c2_2}
\dot{x_2}=-x_2\,(-x_2^2\,y_2^4-y_2^4+a_{32}),\quad 
\dot{y_2}=-y_2^5.
\end{equation}
Then the origin of this system is a stable semi-hyperbolic node (see Theorem 2.15 of \cite{DumLliArt2006}), and going back through the changes of variables we obtain that the origin in the chart $U_2$ is a eliptic sector. See, Figure~\ref{fi:sectoreliptico}. Then there are orbits of the \eqref{center} going to infinity. Hence system \eqref{center} cannot have a global center. 

\begin{figure}[h]
	\begin{center}
		\begin{overpic}[height=4.5cm]{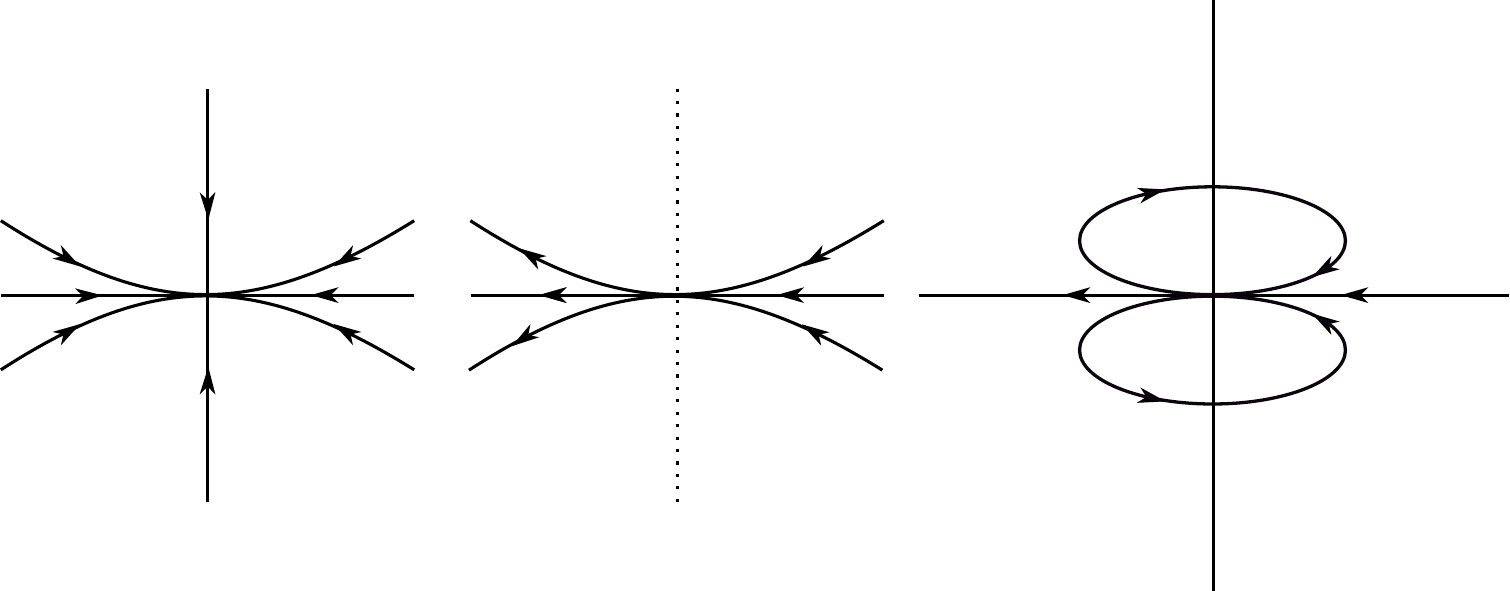}
			\put(25,17){$x_2$}
			\put(14.5,32.5){$y_2$}
			\put(56,17){$x_1$}
			\put(45,32.5){$y_1$}
			\put(98,17){$x$}
			\put(81.5,38){$y$}
		\end{overpic}
	\end{center}
	\caption{In the picture on the left the origin of the system \eqref{c2_2} is a semi-hyperbolic node. In the middle of the figure to the system \eqref{c2_1} all points on the $y_1$-axis are equilibrium points, and by rescaling the time in the left flow change of direction concerning the previous one, and on the right the local phase portrait in the neighborhood of the equilibrim point localized at the origin of the system \eqref{c2}, is a nilpotent elliptic sector.}\label{fi:sectoreliptico}
\end{figure}
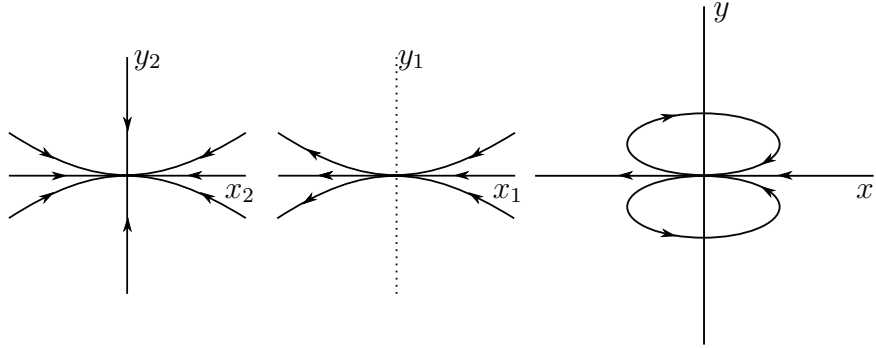

$\bullet$ Case $c_3$. Then there are no infinite singular points in the local chart $U_1$. So it is enough study the origin of the local chart $U_2$. Now system \eqref{U2} is
\begin{equation}\label{c3}
\dot{x}=a^2\,x^6 + x^2\,y^4+ y^4- a_{32}\,x^4,\quad
\dot{y}=x\,y\,(a^2\,x^4+y^4-a_{32}\,x^2).
\end{equation}
Since $\gamma_{4}=y^5=0$, we do the vertical blow-up $(x,y)\rightarrow(x_1,y_1\,x_1),$ and system \eqref{c3} becomes
\begin{equation}\label{c3_1}
\dot{x_1}=x_1^4\,(x_1^2\,y_1^4+y_1^4+a^2\,x_1^2- a_{32}), \quad
\dot{y_1}=-x_1^3\,y_1^5,
\end{equation}
and doing the rescaling $(x_2,y_2,t)\rightarrow(x_1,y_1, t/x_1^3)$ to system \eqref{c3_1} we get
\begin{equation}\label{c3_2}
\dot{x_2}=x_2\,(x_2^2\,y_2^4+y_2^4+a^2\,x_2^2- a_{32}),\quad
\dot{y_2}=-y_2^5.
\end{equation}
Thus the origin of system \eqref{c3_2} is a semi-hyperbolic saddle (by Theorem 2.15 of  \cite{DumLliArt2006}), see the left picture of Figure \ref{fi:sectorhyperbolic}. Going back to system \eqref{c3_1} we obtain the phase portrait of the middle picture of Figure \ref{fi:sectorhyperbolic}, in that picture the $y_1$-axis is filled with equilibria. Undoing the vertical blow-up we obtain the nilpotent hyperbolic sector at the origin of the local chart $U_2$ corresponding to system \eqref{c3} showing in the right picture of Figure \ref{fi:sectorhyperbolic}. Therefore, by Proposition~\ref{cg}, in this case system \eqref{center} has a global center.

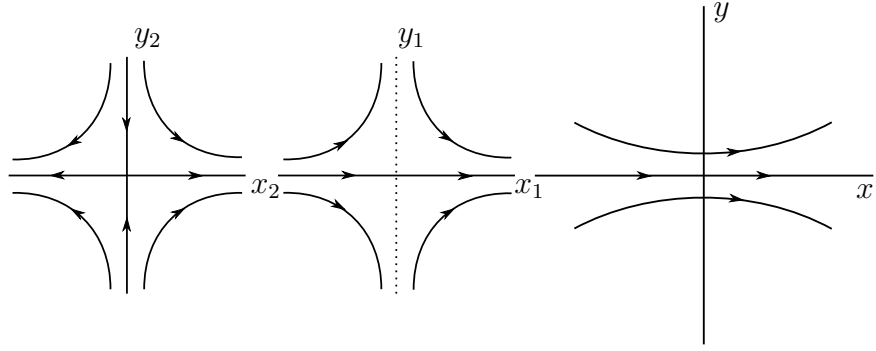
\begin{figure}[h]
\begin{center}
\begin{overpic}[height=4.5cm]{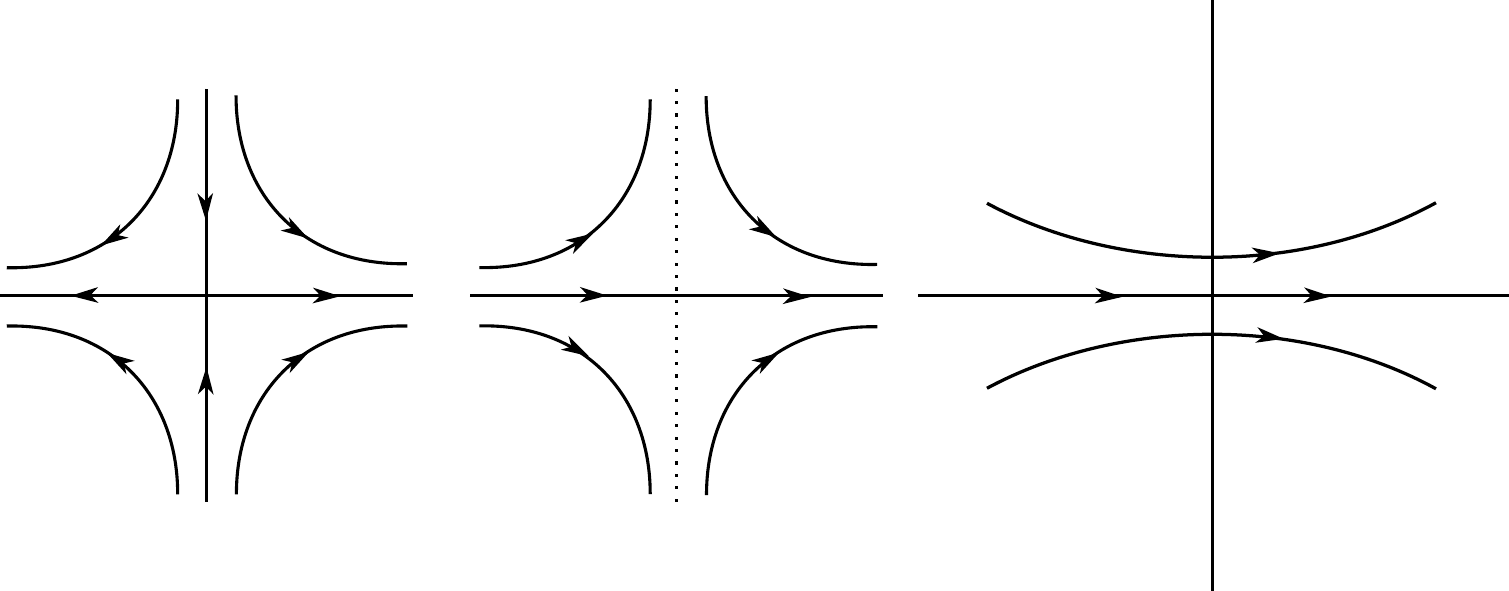}
\put(28,17.5){$x_2$}
\put(14.5,35){$y_2$}
\put(58.5,17.5){$x_1$}
\put(45,35){$y_1$}
\put(98,17){$x$}
\put(81.5,38){$y$}
\end{overpic}
\end{center}
\caption{ The local phase portraits corresponding to the blow-up of the origin in the chart $U_2$ of system \eqref{c3}. %In the left the origin of the systems \eqref{c3_2} is a semi-hyperbolic saddle. In the middle to the system \eqref{c3_1} all points on the $y_2$-axis are equilibrium points, and by rescaling the time in the left flow change of direction concerning the previous one.
}\label{fi:sectorhyperbolic}
\end{figure}

$\bullet$ Case $c_4$. For this case we have the systems \eqref{c3}, \eqref{c3_1} and \eqref{c3_2} with $a=0$, and as in the case $c_3$ the local phase portrait at the origin of the local chart $U_2$ is shown in the right picture of Figure \ref{fi:sectorhyperbolic}. Now we must study the local phase portrait at the origin of the local chart $U_1$. So system \eqref{U1} is
\begin{equation}\label{c4}
\dot{x}=(a_{32}-y^4)\,x^2-y^4,\quad \dot{y}=-y^5\,x.
\end{equation}
Since $\gamma_6=-y\,(y^4-a_{3,2}\,x^2)$ we can do the vertical blow-up $(x,y)\rightarrow(x_1,y_1\,x_1)$ to system \eqref{c4}, and we obtain
\begin{equation}\label{c4_1}
\dot{x_1}=-x_1^2\,(x_1^4\,y_1^4+x_1^2\,y_1^4 -a_{32}),\quad
\dot{y_1}=x_1\,y_1\,(x_1^2\,y_1^4-a_{32}),
\end{equation}
then with the rescaling $(x_2,y_2,t)\rightarrow(x_1,y_1,t/x_1)$, system \eqref{c4_1} writes
\begin{equation}\label{c4_2}
\dot{x_2}=-x_2\,(x_2^4\,y_2^4+x_2^2\,y_2^4 -a_{3,2}),\quad
\dot{y_2}=y_2\,(x_2^2\,y_2^4-a_{3,2}).
\end{equation}
Therefore the origin of system \eqref{c4_2} is a hyperbolic saddle, its phase portrait is the one of the left picture of Figure \ref{fi:sectorhyperbolic} but with the orbit run in reversing sense. Going back through the changes of variables we obtain that the local phase portrait at the origin of the local chart $U_1$ is the one of the right picture of Figure  \ref{fi:sectorhyperbolic}, reversing the orientation of the orbits.  So, again from Proposition~\ref{cg}, in this case system \eqref{center} has a global center.

$\bullet$ Case $c_5$. From \eqref{U1_1} it follows that there are no infinite singular points in the local chart $U_1$. Hence we must study only the origin of the local chart $U_2$. So system \eqref{U2} is
\begin{equation}\label{c5}
\dot{x}=a^2\,x^6 + x^2\,y^4 + y^4,\quad
\dot{y}=x\,y\,(a^2\,x^4 + y^4).
\end{equation}
Since $\gamma_6=y^5$, we do the vertical blow-up and after the rescaling of the time $(x,y,t)\rightarrow(x_1,y_1\,x_1,t/x_1^3),$ and  we obtain
\begin{equation}\label{c5_1}
\dot{x_1}=x_1\,(x_1^2\,y_1^4+y_1^4+a^2\,x_1^2),\quad
\dot{y_1}=-y_1^5.
\end{equation}
Now $\gamma_5=a^2\,x_1^3\,y_1$, therefore $x_1=0,$ is a characteristic direction and before doing a vertical blow-up we do the twist $(x_1,y_1)\rightarrow(x_2-y_2,y_2),$ in system \eqref{c5_1} and we get
\begin{equation}\label{c5_2}
\begin{array}{rl}
\dot{x_2}=&x_2^3\,y_2^4-3\,x_2^2\,y_2^5+3\,x_2\,y_2^6-y_2^7-2\,y_2^5+ x_2\,y_2^4+a^2\,x_2^3-a^2\,(3\,x_2^2\,y_2\\
&-3\,x_2\,y_2^2- y_2^3), \vspace{0.1cm}\\
\dot{y_2}=&-y_2^5.
\end{array}	
\end{equation}
Since $\gamma_5=a^2\,y_2\,(x_2-y_2)^3$, we do the vertical blow-up and the rescaling of the time $(x_2,y_2,t)\to (x_3,y_3\,x_3,t/x_3^2)$ in system \eqref{c5_2}, so we obtain the system
\begin{equation}\label{c5_3}
\begin{array}{rl}
\dot{x_3}=&-x_3\,(x_3^4\,y_3^7-3\,x_3^4\,y_3^6+3\,x_3^4\,y_3^5 - x_3^4\,y_3^4+2\,x_3^2\,y_3^5-x_3^2\,y_3^4+a^2\,y_3^3\\
&-3\,a^2\,y_3^2 + 3\,a^2\,y_3-a^2),\vspace{0.1cm} \\
\dot{y_3}=&(x_3^4\,y_3^6-2\,x_3^4\,y_3^5+x_3^4\,y_3^4+2\,x_3^2\,y_3^4+a^2\,y_3^2-2\,a^2\,y_3+a^2)(y_3^2-y_3).
\end{array}	
\end{equation}

The equilibria of system \eqref{c5_3} on the straight line $x_3=0$, are the origin and the $(0,1)$. The origin is a hyperbolic saddle, and the $(0,1)$ is linearly zero, i.e. the matrix of the linear part of system \eqref{c5_3} is identically zero. Then we will study the local phase portrait at the point $(0,1)$ doing blow-up's. First we translate the equilibrium point $(0,1)$ to the origin of coordinates, so applying in \eqref{c5_3} the change $(x_3,y_3)\rightarrow (x_4,y_4+1),$  we get 
\begin{equation}\label{c5_4}
	\begin{aligned}
	\dot{x_4}=&-x_4(x_4^2+6\,x_4^2+14\,x_4^2\,y_4^2+a^2\,y_4^3+16\,x_4^2\,y_4^3+x_4^4\,y_4^2+9\,x_4^2\,y_4^4+4\,x_4^4\,y_4^4\\
	&+2\,x_4^2\,y_4^5+6\,x_4^4\,y_4^4+4\,x_4^4\,y_4^6+x_4^4\,y_4^7),\\
	\dot{y_4}=&y_4(1+y_4)(2\,x_4^2+8\,x_4^2\,y_4+a^2\,y_4^2+12\,x_4^2\,y_4^2+x_4^4\,y_4^2+8\,x_4^2\,y_4^3+4\,x_4^4\,y_4^3\\&+2\,x_4^2\,y_4^4+6\,x_4^4\,y_4^4+4\,x_4^4\,y_4^5+x_4^4\,y_4^6).
	\end{aligned}
\end{equation}
For system \eqref{c5_4} we have $\gamma_5=-x_4\,y_4(3\,x_4^2+a^2\, y_4^2)$, then $x_4=0,$ is a characteristic direction. Therefore we do the twist $(x_4,y_4)\rightarrow(x_5-y_5,y_5,$ to system \eqref{c5_4}, and after we do the vertical blow-up and the rescaling of the time $(x_5,y_5,t) \rightarrow(x_6,y_6\,x_6,t/x_6^2)$, obtaining the system 
\begin{equation}\label{c6_5}
	\begin{aligned}
		\dot{x_6}=&-x_6\big(1-5\,y_6+f(x_6,y_6,a)\big),\\
		\dot{y_6}=&y_6(1-y_6)\big(3-6\,y_6+g(x_6,y_6,a)\big),
	\end{aligned}
\end{equation}
where $f$ and $g$ are polynomials of degree $21$ and $20$ in the variables $x_6, y_6,$ respectively. 

From \eqref{c6_5} we obtain that the equilibrium points on $u_6=0$ are the origin and $(0,1).$ Both equilibria are hyperbolic sadddles. Therefore going back to system \eqref{c5_3}, we obtain the local phase portrait at the equilibrium $(0,1)$. The steps on the blow down from the differential system 
\eqref{c6_5} until the system \eqref{c5_3} are given in the right column of
Figure \ref{fi:sectorhyperbolic_1}. 

Finally going back to the differential system \eqref{c5}, we get the local phase portraits at the origin of the chart $U_2$. The different steps of this blow down are in the right column of Figure \ref{fi:sectorhyperbolic_1}. Again, in this case from Poposition \ref{cg} the differential system \eqref{center} has a global center.         

\begin{figure}[p]
	\begin{center}
		\begin{overpic}[width=11cm,height=16cm]{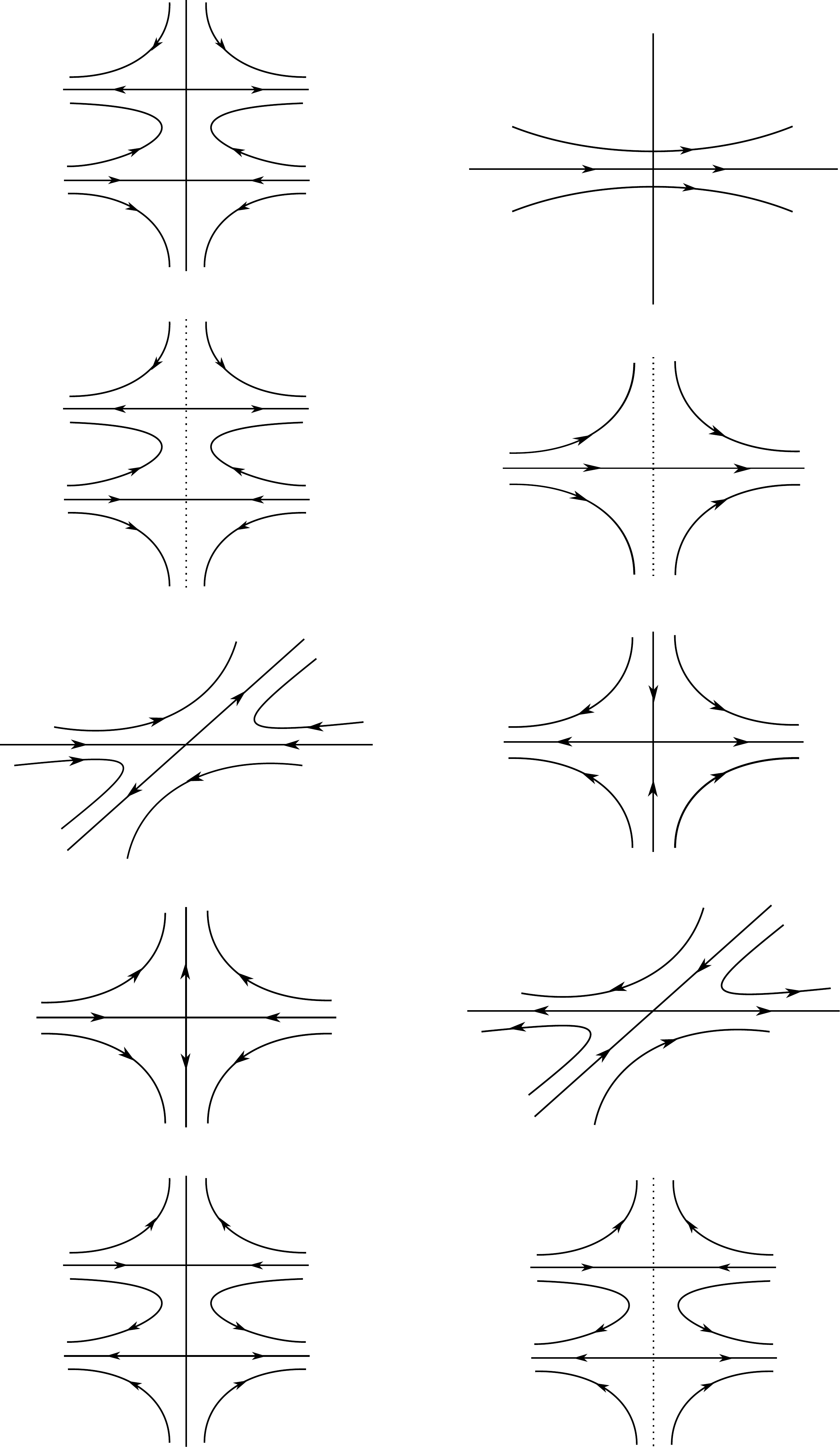}
			\put(26,94){$y_6=1$}
			\put(26,87.5){$x_6$}
			\put(15,101){$y_6$}
			\put(26,72){$y_6=1$}
			\put(26,65.5){$x_6$}
			\put(15,79){$y_6$}
			\put(31,48.5){$x_5$}
			\put(22,56.5){$y_5=x_5$}
			\put(28,30){$x_4$}
			\put(15,38.5){$y_4$}
			\put(25.8,6.5){$x_3$}
			\put(15,20){$y_3$}
			\put(26,13){$y_3=1$}
			\put(64,13){$y_3=1$}
			\put(64,6.5){$x_3$}
			\put(53,20){$y_3$}
			\put(69,30.5){$x_2$}
			\put(60,38){$y_2=x_2$}
			\put(66,49){$x_1$}
			\put(53,57.5){$y_1$}
			\put(66,68){$x_1$}
			\put(53,76.5){$y_1$}
			\put(69,89){$x$}
			\put(53.5,98.5){$y$}
		\end{overpic}	
	\end{center}
\caption{The local phase portrait corresponding to the blow-up of the origin of the chart $U_2$ of system \eqref{c5}. Starting with the top picture on the left column, then going down one by one in the left column we have: The local phase portrait of system \eqref{c6_5}, both equilibrium points the origin and $(0,1)$ are hyperbolic saddles. In the next one all points on $y_6$-axis are equilibrium points. Undoing the blow up we obtain a saddle for system $(\dot x_5,\dot y_5)$, and undoing the twist transformation we continue having a saddle for system $(\dot x_4,\dot y_4)$. }\label{fi:sectorhyperbolic_1}
\end{figure}
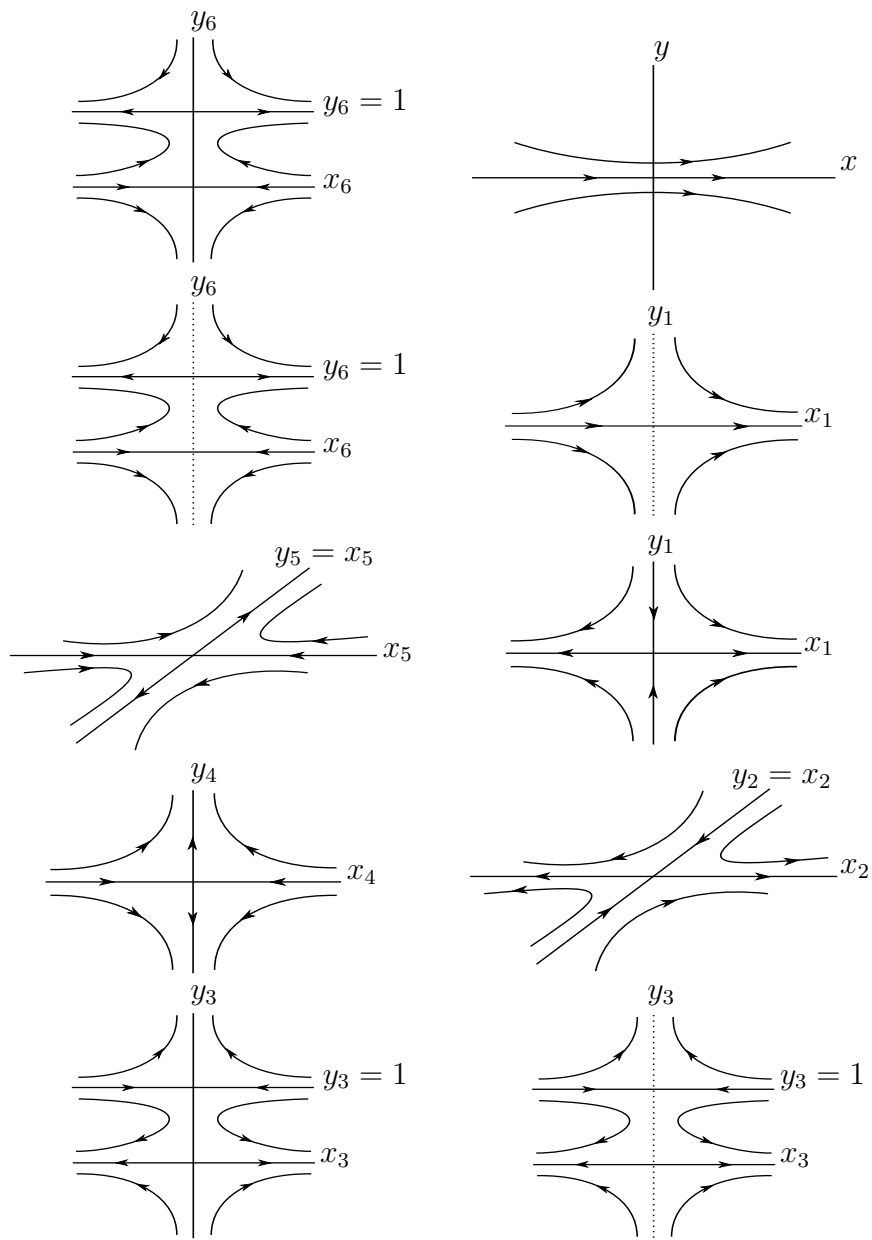

$\bullet$ Case $c_6$. Then the differential system \eqref{eq:1} under the conditions of Theorem \ref{thmmainn_1} becomes the linear differential system $\dot x=y$. $\dot y=-x$.

In summary, from the cases $c_3$, $c_4$ and $c_5,$ it follows the proof of Theorem \ref{thmmainn_2}.
\end{proof}

\noindent\hspace{0.10\linewidth}\begin{minipage}{0.82\linewidth}     In the following we get for system $(\dot x_3,\dot y_3)$ two saddles one hyperbolic at the origin and one linearly zero at $(0,1)$.
	 Now, at the end of the right column, undoing the time rescaling we obtain the same local phase portrait in a neighborhood of $u_3=0$ with the exception that the straight line $u_3=0$ is filled with   equilibrium points. Undoing the blow-up, the twist transformation, the first rescaling and the first blow up we obtain that the origin of $U_2$ is formed by two hyperbolic sectors having their separatrices on the infinite circle.
\end{minipage}

\section*{Acknowledgements}

Leonardo P.C. da Cruz was supported by São Paulo Paulo Research Foundation (FAPESP) grants number 2022/14484-9, 2021/14987-8. The second author is partially suportted by the Agencia Estatal de Investigaci\'on grant PID2019-104658GB-I00, the H2020 European Research Council grant MSCA-RISE-2017-777911, AGAUR (Generalitat de Catalunya) grant 2021SGR00113, and by the Acad\`emia de Ci\`encies i Arts de Barcelona.

\end{document}